\documentclass[11pt] {article} 
  \usepackage{amsmath}
    \usepackage{amssymb}

\newtheorem{example}{Example}[section]
\newtheorem{theorem}{Theorem}[section]
\newtheorem{lemma}{Lemma}[section]

\newtheorem{remark}{Remark}[section]

\newcommand{\eqnsection}{
   \renewcommand{\theequation}{\thesection.\arabic{equation}}
   \makeatletter
   \csname @addtoreset\endcsname{equation}{section} 
   \makeatother}
\def \ov{\overline}

\def \be{\begin{equation}}
\def \ee{\end{equation}}
\def \bt{\begin{theorem}} 
\def \et{\end{theorem}}
\def \bl{\begin{lemma}} 
\def \el{\end{lemma}}
\def \bea{\begin{eqnarray}}
\def \eea{\end{eqnarray}}
\def \bas{\begin{eqnarray*}}
\def \eas{\end{eqnarray*}}

\def \al{\alpha}
\def \bb{\beta}

\def \de{\delta}

\def \ep{\epsilon}

\def \la{\lambda} 
\def \La{\Lambda}

\def \Om{\Omega}

\def \si{\sigma}

\def \th{\theta}

\def \ff{\infty}
\def \wh{\widehat}
\def \wt{\widetilde}

\def \cd{\,\cdot\,}

\def\stl{\stackrel{law}{=}}

\def \DD{{\cal D}}

\def \FF{{\cal F}}

\def \II{{\cal I}}

\def \KK{{\cal K}}

\def \TT{{\cal T}}

\def \VV{{\cal V}}

\def\b1{\mathbf 1}
\def \({\left(}
\def \){\right)}

\def \nn{\nonumber}
\def \Proof{\noindent{\bf Proof $\,$ }}

\def \bc{\begin{center} }
\def \ec{\end{center} }
\def \bs{\begin{slide} }
\def \es{\end{slide} }

\def\square{{\vcenter{\vbox{\hrule height.3pt
        \hbox{\vrule width.3pt height5pt \kern5pt
           \vrule width.3pt}
        \hrule height.3pt}}}}
\def\qed{{\hfill $\square$ \bigskip}}

\eqnsection

\begin{document}

\title{  Permanental processes with kernels that are not equivalent to a symmetric matrix }

  \author{  Michael B. Marcus\,\, \,\, Jay Rosen \thanks{Research of     Jay Rosen was partially supported by  a grant from the Simons Foundation.   }}
\maketitle

 \begin{abstract}       Kernels of $\al$-permanental processes of the form
 \begin{equation}
 \wt  u(x,y)=u(x,y)+f(y),\qquad x,y\in S,
   \end{equation}
 in which $u(x,y)$ is symmetric, and $f$ is an excessive function for the Borel right process with potential densities $ u(x,y)$, are considered. Conditions are given that determine whether $\{\wt   u(x,y);x,y\in S\}$ is symmetrizable or asymptotically symmetrizable.

    \end{abstract}  
 
\maketitle

\section{Introduction}\label{sec-1}

 	 An $R^{n}$ valued   $\al$-permanental random variable $X=(X_{1},\ldots, X_{n})$ is  a random variable with Laplace transform 
\begin{equation}
   E\(e^{-\sum_{i=1}^{n}s_{i}X_{i}}\) 
 = \frac{1}{ |I+KS|^{ \al}},   \label{int.1} 
 \end{equation}
where $K$ is an $n\times n$ matrix and  $S $ is  an   $n\times n$  diagonal matrix with   diagonal entries  $(s_{1},\ldots,s_{n})$.

  We refer to $K$ as a kernel of $X$.   But note that   $K$ is not unique. For example, if $K$ satisfies (\ref{int.1}) so does 
 $\La K\La^{-1}$ for any $\La\in \DD_{n,+}$, the set of  $n\times n$  diagonal matrices with strictly positive diagonal entries.  
 
Let $\KK(X)$ denote  the set of all kernels that determine  $X $ by (\ref{int.1}).
We are  particularly interested in $\al$-permanental random variables $X$ for which $\KK(X)$  does not contain any symmetric kernels.  (We  explain at the end of this section why we are  interested in such processes and kernels.)

If $\KK(X)$ contains a symmetric matrix we say that $X$   is determined by a symmetric matrix or kernel and     that   any $K\subset \KK(X)$ is equivalent to a symmetric matrix, or   is  symmetrizable.    It follows from  (\ref{int.1}) that a kernel $K$ is equivalent to a symmetric matrix if and only if there exists an $n\times n$  symmetric matrix $Q$ such that   
	 \begin{equation}
    |I+KS| = |I+QS| \quad\mbox{for all $S\in \DD_{n,+}$}\label{1.3qqr}.
   \end{equation}

 	An $\al$-permanental process $\{X_{t},t\in T\}$ is a stochastic process that has finite dimensional distributions that are $\al$-permanental random variables.  An $\al$-permanental process is determined by a kernel $\{K(s,t),s,t\in T\}$ with the property that for all distinct $t_{1},\ldots,t_{n}$ in $T$,  $\{K(t_{i},t_{j}),i,j\in [1,n]\}$ is the kernel of the $\al$-permanental random variable $(X_{t_{1}},\ldots,X_{t_{n}})$.

 \medskip \noindent 	{\bf Definition }We say that an $\al$-permanental process $\{X_{t},t\in T\}$ with kernel $\{K(s,t),s,t\in T\}$ is   determined by a symmetric kernel  if  for all $n\ge 1$ and distinct $t_{1},\ldots,t_{n}$ in $T$,  $\{K(t_{i},t_{j}),i,j\in [1,n]\}$ is symmetrizable.  When this is the case we also say that $\{K(s,t),s,t\in T\}$ is symmetrizable. 
	  (In what follows we always take  $|T|\ge 3.)$

 \medskip	The next theorem is \cite[Theorem 1.9]{MRall}. It shows that we can modify a very large class of symmetric potentials  so that they are no longer symmetric   but are still kernels of permanental processes.

 	   \bt
   \label{theo-borelN}
        Let $S$ a be locally compact set with a countable base. 
Let $X\!=\!
(\Om,  \FF_{t}, X_t,\th_{t},P^x
)$ be a transient symmetric Borel right process with state space $S$  and continuous strictly positive  potential   densities  $u(x,y)$ with respect to some $\si$-finite measure $m$ on $S$.  
 Then for any   finite excessive function   $f$ of $X$ and $\al>0$,
\begin{equation}
   \wt u^{f}(x,y)= u(x,y) +f(y),\qquad x,y\in S,\label{1.10mm}
   \end{equation}
is the kernel of an $\al$-permanental  process.
\et

   A function $f$ is said to be  excessive for $X$ if $  E^{x}\(f(X_{t})\)\uparrow  f(x)$ as $t\to 0$ for all  $x\in S$.     
  It is easy to check that for any positive measurable function $h$, \begin{equation}
   f(x)=\int u(x,y) h(y) \,dm(y)=E^{ x} \(   \int_{0}^{ \ff}h\(   X_{t}\)\,dt\)\label{potdef}
   \end{equation}
  is excessive for $X$.   
Such a function $f$   is called  a potential function for $X$.   

\medskip	 Unless the   function   $f$ in  (\ref{1.10mm}) is constant, 
 $\{\wt u^{f}(x,y);x,y\in S\}$ is not symmetric. We now show that, generally, we can choose $f $ so that  $\{\wt u^{f}(x,y);x,y\in S\}$ is also not equivalent to a symmetric matrix.     The next two theorems show how restricted the symmetric matrix $\{  u(x,y);x,y\in S\}$ must be for $\{\wt u^{f}(x,y);x,y\in S\}$ to be symmetrizable  for all potential functions $f$.  

\medskip	We  use $\ell_{1}^{+}$ to denote strictly positive sequences   in $\ell_{1}$.   

 \bt\label{theo-borelNS}
Let $X\!=\!
(\Om,  \FF_{t}, X_t,\th_{t},P^x
)$ be a transient symmetric Borel right process with state space   $  T\subseteq \mathbb N$, and potential  $U=\{U_{j,k}\}_{j,k\in   T}$. Then  
\begin{itemize}
\item [(i)] Either
 \begin{equation}
 U_{j,k} =\La_{j}\de_{j,k}+ d, \qquad j,k\in    T,\label{1.5nn}
   \end{equation}
   where $\La_{j}\ge 0$ and $d\ge 0$,
   \item [(ii)] or we can find a potential function   $f=Uh$, with $h  \in \ell^{+}_{1}$,  such that 
   \begin{equation}
  \wt U_{j,k}^{f}:= U_{j,k}+ f_{k},\qquad j,k\in    T,\label{1.6nn}
   \end{equation}
 is not symmetrizable.
    \end{itemize}
%{\bf We do not mention  this  in the proofs. Can we just leave it out? }This dichotomy also  holds  %if   the potential $f=Uh$ is replaced by $f+q$, for any constant $q>0$.
 \et

 When we consider limit theorems for infinite sequences of permanental random variables $\{Y(k), k\in \mathbb N\}$ with kernel $V=\{v(j,k), j,k\in \mathbb N\}$ it is not enough to know that $V$ is not symmetrizable since we are only concerned with 
the permanental variables generated by $V(n)=\{v(j,k), j,k\ge n  \}$ as $n\to \ff$. We would like to know that $V(n)$ is not symmetrizable for   large $n$. We   say that   the kernel $V$ is asymptoticly symmetrizable  if  there exists an $n_{0}$ such that $V(n)$ is symmetrizable for all $n\ge n_{0}$. We can modify  Theorem \ref{theo-borelNS} to handle this case also.     

\bt\label{theo-borelNSmm}
Let $X\!=\!
(\Om,  \FF_{t}, X_t,\th_{t},P^x
)$ be a transient symmetric Borel right process with state space $ \mathbb N$, and potential  $U=\{U_{j,k}\}_{j,k\in   \mathbb N}$. Then  
\begin{itemize}
\item [(i)] Either there exists an $n_{0}$ such that 
 \begin{equation}
 U_{j,k} =\La_{j}\de_{j,k}+ d,\qquad \forall j,k\ge n_{0},\label{3.3mma}
   \end{equation}
   where $\La_{j}\ge 0$ and $d\ge 0$,  
   \item [(ii)] or we can find a potential function   $f=Uh$, with $h  \in \ell^{+}_{1}$,  such that 
   \begin{equation}
   \wt U_{j,k}^{f}:= U_{j,k}+  f_{k},\qquad j,k\in  \mathbb N, 
   \end{equation}
   is not asymptoticly symmetrizable.  
   \end{itemize}
%{\bf As above !!! }This dichotomy also holds  if   the potential $f=Uh$ is replaced by $f+q$, for any %constant $q>0$.
 \et
 
  The next theorem shows that when the state space of a transient symmetric Borel right process has a limit point, then under reasonable conditions on the potential densities that determine the process, the process is not    determined by a kernel that is asymptoticly symmetrizable.

\begin{theorem}\label{theo-1.4}  Let   $S'=\{x_{0},x_{1},\ldots\}$ be a countable set with a single limit point $x_{0}$.  Let $\ov X$ be a transient symmetric Borel right process with state space $ S'$, and continuous strictly positive  potential   densities   $u:=\{u(x,y), x,y\in S'\}$ such that $u(y, x_{0})<u(x_{0},x_{0})$ for all $y\neq x_{0}$.        Then  we can find a potential function     $f=uh$, with  $h  \in \ell ^{+}_{1}$,   that is continuous at $x_{0}$, and is such that,    
 \begin{equation}
   \wt u^{f}(x,y)= u(x,y)+f(y),\qquad x,y\in S',\label{1.9mm}
   \end{equation}
   is not   asymptoticly symmetrizable.     
   
%{\bf As above !!! }   This also holds  if   the potential $f=Uh$ is replaced by $f+q$, for any constant %$q>0$.
 \et
 
  	   Theorems \ref{theo-borelNS}--\ref{theo-1.4} show that generally there exists an excessive function $f$ for $X$ which gives a kernel for an $\al$-permanental processes that is not determined by a   symmetric matrix. However, in specific examples we deal with    specific functions $f$ and want to know that the kernels determined by these functions are not symmetrizable. With some additional structure on the symmetric matrix $u(x,y)$ in (\ref{1.10mm}) we can show that $\wt u^{f}(x,y)$ in (\ref{1.10mm}) is not  asymptoticly symmetrizable.

\begin{lemma}\label{lem-1.1mm} In the notation of (\ref{1.10mm}),   let   $u=\{u(j,k); j,k\in \mathbb N \}$ be a symmetric  To\"eplitz matrix,   with at least two different off diagonal elements, and set $v(|j-k|)=u(j,k)$. Let\begin{itemize}
\item[(i)]\begin{equation}
   \wt u^{f}(j,k)= v(|j-k|) +f(k),\qquad j,k\in \mathbb N ,\label{1.10a}
   \end{equation}
where    $f$ is a strictly monotone potential for $u$. Then $\{\wt u^{f}(j,k);j,k\in \mathbb N\}$ is not asymptoticly symmetrizable.
\item[(ii)]Let
 \begin{equation}
   \wt v^{f}(s_{j},s_{k})= s_{j}\wedge s_{k} +f(s_{k}),\qquad  {j}, {k}\in \mathbb N,\label{1.10b}
   \end{equation}
where   $f$ is a strictly monotone potential for $\{s_{j}\wedge s_{k}; {j}, {k}\in \mathbb N\}$.   Then  for any   triple of  distinct values $ s_{j}, s_{k},s_{l} $, 
    \begin{equation}
  \{ \wt v^{f}(s_{p},s_{q})\}_{p,q=j,k,l}\, ,\label{1.10bb}
   \end{equation}
  is   not  symmetrizable.     In particular  $\{  \wt v^{f}(s_{j},s_{k};j,k\in \mathbb N\}$ is   not  asymptoticly symmetrizable. \end{itemize}
 \end{lemma}

  We can use this lemma to  show that certain $\al$-permanental processes,  studied in \cite{MRall},   are not   determined by kernels that are asymptoticly symmetrizable.    When $S$ is an interval on the real line we say that   $ \{u(x,y);x,y\in S\}$ is not    asymptoticly symmetrizable at  $x_{0}\in S$,  if we can find a sequence $\{x_{k}\}$ in $S$ such that $\lim_{k\to\ff}x_{k}=x_{0}$, and 
 $ \{u(x_{j},x_{k});j,k\in \mathbb N\}$ is not asymptoticly symmetrizable.

 \begin{example}\label{ex-1.1} {\rm 
 In \cite[Example 1.3]{MRall} we obtain a limit theorem for the asymptotic behavior of  the sample paths at 0 of $\al$-permanental   processes with the   kernel,  
  \be
  \wh u^{f}(s,t)=e^{-\la |s-t|}+f(t),\qquad s,t\in [0,1],\label{1.33j}
  \ee  
where $f=q+t^{\bb}$, $\bb>2$, and $q\ge q_{0}(\bb)$, a constant depending on $\bb$. We show in Section \ref{sec-4} that   $  \wh u^{f}(s,t)$ is not   asymptoticly symmetrizable at any $s_{0}\in S$.

Similarly  
 \be
  \ov u^{f}(j,k)=e^{-\la |j-k|}+f(k),\qquad j,k\in \mathbb N ,\label{1.33k}
  \ee  
is not    asymptoticly symmetrizable.  
 }\end{example}

\begin{example}\label{ex-1.2} {\rm 
In \cite[Example 1.4]{MRall} we obtain   limit theorems for the asymptotic behavior of  the sample paths at zero and infinity  of $\al$-permanental   processes with the   kernel, 
\begin{equation}
 \wt v^{f}(s,t)=s  \wedge t+f(t),\qquad s,t\ge 0,\label{1.39j}
  \end{equation}
where $f$ is a concave strictly increasing function.  We show in Section \ref{sec-4} that  for any $s_{0}\in R^{+}$ and any sequence of distinct values $\{s_{k}\} $ such that $\lim_{k\to\ff}s_{k}=s_{0}$, $  \wt v^{f}(s_{j},s_{k})$ is   not  asymptoticly symmetrizable.  

  In addition,  
\begin{equation}
 \ov v^{f}(j,k)=j  \wedge k+f(k),\qquad j,k\in \mathbb N ,\label{1.39k}
  \end{equation}
is    not   asymptoticly symmetrizable.
}\end{example} 
 
  We explain why we are particularly interested in $\al$-permanental processes determined by kernels $K$ that are not equivalent to a symmetric matrix.  When $\{u(s,t);s,t\in \TT\}$ is symmetric and is a kernel that  determines   $\al$-permanental processes,  $Y_{\al}=\{Y_{\al}(t),t\in \TT\}$,   then 
 \begin{equation}
   Y_{1/2} \stl \{G^{2}(t)/2,t\in \TT\},
   \end{equation}
  where $G=\{G  (t) ,t\in \TT\}$ is a mean zero Gaussian process with covariance $u(s,t)$. 
  
   If $\al=m/n$ for integers $m$ and $n$,
   \begin{equation}
   Y_{m/n} \stl \sum_{j=1}^{m}\sum_{k=1}^{n}Y_{1/(2n)}^{(j,k)},
   \end{equation}
 where $Y_{1/(2n)}^{(j,k)}$ are independent copies of  $Y_{1/(2n)} $. Therefore, in some sense, $Y_{m/n}$, is only a modification of  the Gaussian process $G$. 
This is not true when  the kernel of  $\al$-permanental processes is not symmetrizable. In this case we get a new class of processes. These are the processes that we find   particularly interesting.

  \medskip		 To study  permanental processes with kernels that are not equivalent to a symmetric matrix our first step is to characterize those kernels that are  equivalent to a symmetric matrix. This is done in Section \ref{sec-2}. In Section \ref{sec-3} we give the proofs  of Theorems \ref{theo-borelNS}--\ref{theo-1.4}. In Section \ref{sec-4} we give the proof of Lemma \ref{lem-1.1mm} and details about Examples \ref{ex-1.1} and \ref{ex-1.2}.

   \section{Kernels that are equivalent to a symmetric matrix }\label{sec-2}

Let $M$ be an $n\times n$ matrix. For $\II\subseteq [1,\ldots,n]$ we define  $M_{\II}$ to be   the $|\II | \times |\II |$ matrix  $\{M_{p,q}\}_{ p,q\in \II}$.   (Recall that $\DD_{n,+}$ is the set of all $n\times n$ diagonal matrices with strictly positive diagonal elements.) 

\begin{lemma}\label{lem-1.1n} Let $K$ be an $n\times n$ matrix  and assume that   
\begin{equation}
    |I+KS| = |I+QS| \quad\mbox{for all  $S\in \DD_{n,+}$}. \label{1.3qq}
   \end{equation}
   Then for all $\II\subseteq [1,\ldots,n]$
   \begin{equation}
    |K _{\II}  | =    |Q _{\II}  |.\label{1.3qq5}
   \end{equation}
In particular
   \begin{equation}
     |K  | =    |Q   |  \label{1.3qq7}
   \end{equation}
   and
  \begin{equation}
   K_{j, j} =Q_{j,j} \quad\mbox{for all}\quad j=1,\ldots n.\label{1.3qq2}
  \end{equation} 
Furthermore, if $Q$ is symmetric, then 
\begin{equation}
   |Q_{j,k}|=(K_{j,k}K_{k,j})^{1/2}\quad\mbox{for all}\quad  i,j =1,\ldots,n  \label{1.5w}
   \end{equation}
   and  
for all distinct $i_{1},i_{2},i_{3}\in [1,\ldots,n]$  
   \begin{equation}
   K_{i_{1},i_{2}}  K_{i_{2},i_{3}}  K_{i_{3},i_{1}}  =K_{i_{1},i_{3}}   K_{i_{2},i_{1}} K_{i_{3},i_{2}}.\label{1.7w}
     \end{equation}
 \end{lemma}

\Proof  Denote the diagonal elements of $S$  by $\{s_{i}\}_{i=1}^{n} $. Let $s_{i}\to 0$ for all $s_{i}\in \II^{c}$ in (\ref{1.3qq})  to get
\begin{equation}
    |I+K_{\II}S| = |I+Q_{\II}S| \quad\mbox{for all $S\in \DD_{|\II|,+}$}. \label{1.3qq1}
   \end{equation}
   Multiply  both sides of (\ref{1.3qq1}) by $|S^{-1}|$ and 
  let  the diagonal components of $S$ go to infinity to get (\ref{1.3qq5}). The relationships in (\ref{1.3qq7}) and (\ref{1.3qq2}) are simply examples of 
(\ref{1.3qq5}).

Let $\II=\{j,k\}$. It follows from (\ref{1.3qq5})   that 
\begin{equation}
K_{i,i}K_{j,j}-K_{i,j}K_{j,i}=Q_{i,i}Q_{j,j}-Q^{2}_{i,j},\label{1.3qq6}
\end{equation}
which by (\ref{1.3qq2}) implies that $K_{i,j}K_{j,i}=Q^{2}_{i,j}$. This  gives (\ref{1.5w}).

Finally, let  $\II=\{i_{1},i_{2},i_{3}\}$ and take the determinants $     |K (\II)|  $ and $ |Q (\II)  |$.  
 It follows from  (\ref{1.3qq5}), (\ref{1.3qq2}) and  (\ref{1.5w}) that  
\bea
&& K_{i_{1},i_{2}}  K_{i_{2},i_{3}}  K_{i_{3},i_{1}} +K_{i_{1},i_{3}}    K_{i_{2},i_{1}}K_{i_{3},i_{2}}\nn \\
&&\hspace{1in}
   =Q_{i_{1},i_{2}}  Q_{i_{2},i_{3}}  Q_{i_{3},i_{1}} +Q_{i_{1},i_{3}}  Q_{i_{2},i_{1}}  Q_{i_{3},i_{2}}\nn\\
  &&\hspace{1in} =2Q_{i_{1},i_{2}}  Q_{i_{2},i_{3}}  Q_{i_{3},i_{1}}.
\eea
By (\ref{1.5w}) this is equal to 
\begin{equation}
\pm 2(K_{i_{1},i_{2}}  K_{i_{2},i_{3}}  K_{i_{3},i_{1}} K_{i_{1},i_{3}} K_{i_{2},i_{1}}K_{i_{3},i_{2}})^{1/2}   .\label{}
\end{equation}
Set  \begin{equation}
x=K_{i_{1},i_{2}}  K_{i_{2},i_{3}}  K_{i_{3},i_{1}}\quad\mbox{and}\quad    y=K_{i_{1},i_{3}}    K_{i_{2},i_{1}}K_{i_{3},i_{2}}. \label{}
\end{equation}
Then we have
\begin{equation}
x+y=\pm 2\sqrt{xy}. \label{234}
\end{equation}
It is clear from this that $x$ and $y$ have the same sign. If they are both positive, we have
\begin{equation}
x+y= 2\sqrt{xy}, \label{235}
\end{equation}
That is, $(\sqrt{x}-\sqrt{y})^{2}=0$, which gives (\ref{1.7w}).

On the other hand, if $x $ and $y$ are both negative, (\ref{234}) implies that
\begin{equation}
(-x)+(-y)= 2\sqrt{(-x)(-y)}, \label{236}
\end{equation}
which also gives (\ref{1.7w}).\qed

% \begin{remark} \label{rem-2.1}{\rm  \nc The results in Lemma \ref{lem-1.1n} give similar criteria for $K^{-1}$. To see this 
%note that if $K$ is invertible, then by  (\ref{1.3qq7}), $Q$  is also invertible. Set  $A=K^{-1}$, and multiply  (\ref{1.3qq}) by $|A|=|K|^{-1}=|Q|^{-1}$. This gives
%\begin{equation}
%    |S+A| = |S+Q^{-1} | \quad\mbox {for all  $S\in\DD_{n,+}$}.\label{1.3ee}
%   \end{equation} Next multiply these determinants by  $|S^{-1}|$ to get  
%   \begin{equation}
%    |I+AS^{-1}| = |I+Q^{-1}S^{-1}| \quad\mbox{for all  $S\in\DD_{n,+}$}.\label{1.3ee1}
%   \end{equation}
%Therefore, all the results of Lemma \ref{lem-1.1n} hold with $K$ and $Q$ replaced by $A $ and $Q^{-1}$. In particular for $K$ to be equivalent to a symmetric matrix it is necessary that for all distinct $i_{1},i_{2},i_{3}\in [1,\ldots,n]$, 
%    \begin{equation}
%   A_{i_{1},i_{2}}  A_{i_{2},i_{3}}  A_{i_{3},i_{1}}  =A_{i_{1},i_{3}}    A_{i_{3},i_{2}}A_{i_{2},i_{1}}.\label{1.7ww}
%     \end{equation}  
%    }\end{remark}
% 
% 

\begin{remark} {\rm Even when $K$ is  the kernel of $\al$-permanental processes we must have absolute values on the left-hand sides of (\ref{1.5w}). This is because when (\ref{1.3qq}) holds it also holds when $|I+ QS |$ is replaced by $|I+\VV Q\VV S|$ for any signature matrix $\VV$. (A signature matrix is a diagonal matrix with diagonal entries $\pm 1$.) So the symmetric matrix $Q$ need not be the kernel of $\al$-permanental processes   On the other hand, by \cite[Lemma 4.2]{EK}, we can find a symmetric matrix $\wt Q$ that is the kernel of $\al$-permanental processes  such that (\ref{1.3qq}) holds with $Q$ replaced by $\wt Q$ and we have $ \wt Q_{j,k} = (K_{j,k}K_{k,j})^{1/2}$.  }\end{remark}

 	\section{Proofs  of Theorems \ref{theo-borelNS}--\ref{theo-1.4}}\label{sec-3}

We begin with a simple observation that lies at the heart of the    proofs of  Theorems \ref{theo-borelNS} and \ref{theo-borelNSmm}. 

 For $y\in R^{n}$ we use $B_{\de}(y)$ to denote a Euclidean ball of radius $\de$ centered at $x$.

\begin{lemma}\label{lem-3.1mm} Let $W=\{w _{j,k}; j,k=1,2,3\}$ be a  positive symmetric matrix  such that   $w _{j,k}\leq  w_{j,j}\wedge w _{k,k}$.  For any  $x=(x_{1},x_{2},x_{3})$ let $\wt W^{x}$ be a $3\times 3$ matrix defined by 
\begin{equation}
   \wt W^{x}_{j,k}=w _{j,k}+x_{k},\qquad j,k=1,2,3.
   \end{equation}
Suppose that  $\wt W^{x}$ is symmetrizable for all  $x\in B_{\de}(x_{0})$, for some $x_{0}\in R^{3}$ and $\de>0$. Then, necessarily,  
\begin{equation}
   w _{j,k}=\La_{j}\de_{j,k}+ d, \qquad j,k=1,2,3,\label{3.2mm}
   \end{equation}
   where $\La_{j}\ge 0$ and $d\ge 0$.
 \end{lemma}

\Proof  It follows from Lemma \ref{lem-1.1n} that for all $x\in B_{\de}(x_{0})$  
\be   \( w_{1,2}+x_{2}\) \( w_{2,3}+x_{3}\) \( w_{3,1}+x_{1}\)  = \( w_{1,3}+x_{3}\) \( w_{2,1}+x_{1}\)
\( w_{3,2}+x_{2}\).  \label{3.3mm}
   \ee  
   We differentiate each side of  (\ref{3.3mm}) with respect to $x_{1}$ and $x_{2}$ in $B_{\de}(x_{0})$ and see that
   \begin{equation}
   w_{2,3}+{x}_{3} =w_{1,3}+x_{3}.
   \end{equation}
Therefore, we must have $  w_{2,3}  =w_{1,3}$.
 Differentiating twice more with respect to $x_{1}$ and $x_{3}$, and $x_{2}$ and $x_{3}$, we
  see that if  (\ref{3.3mm}) holds for all $x\in B_{\de}\(  x_{0}\)$  then
   \begin{equation}
 w_{2 ,3}  = w_{1 ,3},\quad   w_{1,2}  = w_{3,2} ,\quad\mbox{and}\quad  w_{3,1}  = w_{2,1}.   \label{}
   \end{equation}
 This implies that   for some $(d_{1},d_{2},d_{3})$
   \be  W= \left (
\begin{array}{ c cccc }  
w_{1,1}   &d_{2}& d_{3 }   \\
 d_{1} & w_{2,2} & d_{3}   \\
  d_{1} &d_{2}& w_{3,3}  \end{array}\right ).\label{nsz.4}
      \ee 
Furthermore, since $W$ is symmetric, we must have $d_{1}=d_{2}=d_{3}$.
 
  Set $d=d_{i}$, $i=1,2,3$. Then, since $w_{i,i}\ge w_{i,j}$, $i,j=1,2,3$,        we can write $w_{i,i}=\la_{i}+d$ for some $\la_{i}\geq 0$,  $i=1,2,3$.    This shows that (\ref{3.2mm}) holds.\qed

 In using  Lemma \ref{lem-3.1mm} we often consider $3\times 3$ principle submatrices of a larger matrix. Consider the matrix $\{W(x,y)\}_{x,y\in S}$, for some index set $S$. Let $\{x_{1},x_{2},x_{3}\}\subset S$.   Consistent with the notation introduced at the beginning of Section \ref{sec-2} we note that
\begin{equation}
   W_{\{x_{1},x_{2},x_{3}\}}=\{ W_{x_{j},x_{k} }\}_{j,k=1}^{3}.
   \end{equation}
   We also use $1_{n}$ to denote an $n\times n$ matrix with all its elements equal to 1.

  \medskip	 \noindent{\bf  Proof of  Theorem  \ref{theo-borelNS}.} 
 If $(i)$ holds then  
  \begin{equation}
   \wt U^{f}:= \La+1_{|T|}G,
   \end{equation}
where  $G$ is a $|T|\times |T|$  diagonal matrix with   entries $f_{1}+d,f_{2}+d,\ldots  $.  Let 
$\II$ be any finite subset of $T$. Obviously, 
  \begin{equation}
   \(\wt U^{f}\)_{\II}= \La_{\II}+1_{|\II|}G_{\II}.
   \end{equation}
Since  
 \begin{equation}
 G_{\II}^{1/2} \(\La_{\II}+1_{|\II|}G_{\II}\)G_{\II}^{-1/2}= \La_{\II}+G_{\II}^{1/2}1_{|\II|}G_{\II}^{1/2}, 
   \end{equation}
and $\La_{\II}+G^{1/2}_{\II}1_{|\II|}G^{1/2}_{\II}$ is symmetric, we see that $  \wt U^{f}$ is   symmetrizable. This shows that if $(i)$ holds then $(ii)$ does not hold.

Suppose that $(i)$ does not hold.     We    show that in this case
we can find   a triple  $\{t_{1},   t_{2},  t_{3}\}$ such that   $U_{\{t_{1},   t_{2},  t_{3}\}}$
 does not have all its off diagonal elements equal.  

Since $(i)$ does not hold
 there are two off diagonal elements of $V$ that are not equal, say   $u_{l,m}=a$
 and $u_{p,q}=b$. Suppose that none of the indices $l,m,p,q$ are equal. The kernel of $(X_{l},X_{m},X_{p})$ has the form.
 \be U_{\{l,m,p\}} =    \left (
\begin{array}{ c cccc   }  
\cd &a& \cd      \\
a & \cd&\cd        \\
\cd &\cd& \cd    
   \end{array}\right ),\label{3.26x}
      \ee
where we use $\cd$ when we don't know the value of the entry. If any of   the off diagonal terms of $U_{\{l,m,p\}}$ are not equal   to $a$  we are done.

  Assume then that  all the off diagonal terms of  $U_{\{l,m,p\}}$ are equal. This implies, in particular, that $(U_{\{l,m,p\}})_{m,p}=(U_{\{l,m,p\}})_{p,m}=a$. Therefore,    $U_{\{m,p,q\}}$ has the form,
 \be U_{\{m,p,q\}}:=    \left (
\begin{array}{ c cccc   }  
\cd &a& \cd      \\
a & \cd&b       \\
\cd &b& \cd    
   \end{array}\right ).\label{3.26}
      \ee
       Therefore, if none of the indices $l,m,p,q$ are equal we see that there exists a triple     $\{t_{1},   t_{2},  t_{3}\}$ such that $U_{\{t_{1},   t_{2},  t_{3}\}}$ does not have all its off diagonal elements equal.

If $l=p$   the argument  is simpler,  because   in this case  
\be    U_{\{l,m,q\}} = \left (
\begin{array}{ c cccc   }  
\cd &a& b      \\
a & \cd&\cd       \\
b &\cd & \cd    
   \end{array}\right ).\label{3.27}
      \ee
  If $m=q$ the kernel of  $(X_{l},X_{p}, X_{m})$ is
\begin{equation}\(
   \begin{array}{ c cccc   }  
\cd &\cd& a     \\
\cd & \cd&b      \\
a &b & \cd    
   \end{array}\right ).\label{3.2ww}
      \ee
  Using the fact that $U$ is symmetric we see that cases when   $l=q$ or $m=p$ are included in the above.   

  This  shows that when $(i)$ does not hold we can find  a triple  $\{t_{1},   t_{2},  t_{3}\}$ such that $U_{\{t_{1},   t_{2},  t_{3}\}}$ does not have all its off diagonal elements equal. We now show that in this case $(ii)$ holds, that is, we can find a potential $f$ for which (\ref{1.6nn}) is not symmetrizable.

 			For convenience we rearrange the   indices    so that   $\{t_{1},   t_{2},  t_{3}\}$ =  $\{1,  2,  3\}$.  We take any $h^{*}\in \ell_{1}^{+}$ and consider the potential $f^{*}=U h^{*}$.  If   $U_{1,2,3}:= \{U_{j,k}+f^{*}_{k}\}_{j,k=1}^{3}$ is not symmetrizable, we are done.   That is, $(ii)$ holds with $f=f^{*}$. However, it is possible that   $U_{\{1,2,3\}}$ is not of the form  of  (\ref{3.2mm}) but   
\be   \( U_{1,2}+f^{*}_{2}\) \( U_{2,3}+f^{*}_{3}\) \( U_{3,1}+f^{*}_{1}\)  = \( U_{1,3}+f^{*}_{3}\) 
\( U_{2,1}+f^{*}_{1}\)\( U_{3,2}+f^{*}_{2}\). \label{3.3k}
   \ee  
(See (\ref{3.3mm})). Nevertheless, since  $U_{\{1,2,3\}}$ is not of the form  (\ref{3.2mm}), it follows from Lemma \ref{lem-3.1mm} that for all $\de>0$
there exists an $(f_{1},f_{2},f_{3})\in  B_{\de}(f^{*}_{1},f^{*}_{2},f^{*}_{3})$ such that $\{U_{j,k}+f  _{k}\}_{j,k=1}^{3}$ is not symmetrizable.
(Here we use  the facts that   a symmetric potential density $U_{j,k}$ is always positive and satisfies $U_{j,k} 	\leq U_{j,j}\wedge U_{k,k}$, see \cite[(13.2)]{book}.)

  Note that   $U_{\{1,2,3\}}$ is invertible. (See e.g.,  \cite[Lemma A.1]{MRnec}.)  Therefore, we can find   $c_{1},c_{2},c_{3}$ such that 
\begin{equation}
   f_{j}=  f_{j}^{*}+ \sum_{k=1}^{3} U_{j,k} c_{k},\qquad j=1,2,3. \label{3.14mm}
   \end{equation}
Now, set  $h=h^{*}+c$, where $c=(c_{1},c_{2},c_{3},0,0,\ldots)$, i.e., all the components of $c$ except for the first three are equal to 0 and set  $f=Uh$.  The components $f_{1},f_{2},f_{3}$ are given by (\ref{3.14mm}). Furthermore, we can choose $\de$ sufficiently small so that for $(f_{1},f_{2},f_{3})\in  B_{\de}(f^{*}_{1},f^{*}_{2},f^{*}_{3})$, $c_{1},c_{2},c_{3}$ are   small enough so that  $h_{1}$, $h_{2}$, and $h_{3}$ are strictly greater than 0, which, of course, implies that $h\in \ell_{1}^{+}$, (defined just prior to Theorem \ref{theo-borelNS}). Therefore, $(ii)$ holds with this potential $f$. \qed

  In Theorem \ref{theo-borelNS} it  is obvious that if $(i)$ does not hold then there are  functions $f$ for which   (\ref{1.6nn}) is not symmetrizable. What was a little difficult was to show that $f=(f_{1},f_{2},\ldots)$,   is a potential for $X$. We have the same problem in the proof of Theorem \ref{theo-borelNSmm} but it is much more complicated. If we start with a potential $f^{*}=Uh^{*}$,  to show that $\wt U^{f}$   is not asymptotically symmetrizable, we may need to modify an infinite number of  the components of $f^{*}$ and still  end up with a potential $f$. The next lemma is the      key to doing this.
 
\bl\label{lem-borelNSmm}
Let $X\!=\!
(\Om,  \FF_{t}, X_t,\th_{t},P^x
)$ be a transient symmetric Borel right process with state space $ \mathbb N$, and potential  $U=\{U_{j,k}\}_{j,k\in \mathbb N}$. Then  we can find a potential function   $f=Uh$, with $h  \in \ell^{+}_{1}$,  such that for all  $\al>0$,
\begin{equation}
   \wt U^{f}_{j,k}= U_{j,k} +f_{ k},\qquad j,k\in\mathbb N,\label{1.10}
   \end{equation}
is the kernel of an $\al$-permanental sequence.

Moreover,  for $I_{l}=\{3l+1,3l+2,3l+3\}$,  the following dichotomy holds for each $l\ge 0$:\begin{itemize}
\item [(i)] Either    $\wt U^{f}_{I_l}$ is not  symmetrizable, 
\item [(ii)]  or   
 \begin{equation}
 U_{I_l}=\La+d  1_{3},\label{nsz.00}
   \end{equation} 
   where  $\La \in D_{3,+}$ and $d\geq 0$.

\end{itemize} 
 \el
 
\Proof  Let $\{i_{l,j}=3l+j\}_{l\ge 0,j\in \{1.2.3\}}$. For $f=\{f_{k}\}_{k=1}^{\ff}$ define,  
  \bea
F_{l}( f )&=&F_{l}( f_{i_{l,1}}, f_{i_{l,2}}, f_{i_{l3}})\label{nsz.2}\\
&=& ( U_{i_{l,1},i_{l,2}}+f_{i_{l,2}})( U_{i_{l,2}, i_{l,3}} +f_{i_{l,3}})( U_{i_{l,3},i_{l,1}}+ f_{i_{l,1}}) \nn\\
 &&\hspace{.1in} -( U_{i_{l,1},i_{l,3}}+f_{i_{l,3}})( U_{i_{l,3},i_{l,2}}+f_{i_{l,2}})(U_{i_{l,2}, i_{l,1}}+ f_{i_{l,1}}). \nn
   \eea 
 We note that when  $U_{I_l}$ is given by (\ref{nsz.00}), then for any sequence  $\{f_{i_{1}}, f_{i_{2}}, f_{i_{3}}\}$, $ F_{l}( f)=0$ and $\wt U^{f}_{I_l}$ is   symmetrizable. The first assertion in the previous sentence follows   because
 all the terms   $\{U_{i_{j}.i_{k}}\}_{j\ne k=1}^{3}$ are equal $d$.    The second 
 is proved in the first paragraph of the proof of  Theorem  \ref{theo-borelNS}.   On the other hand, it follows from Lemma \ref{lem-1.1n} that if $ F_{l}( f)\neq 0$ then $\wt U^{f}_{I_l}$ is not symmetrizable.
      
Therefore, to prove this theorem it suffices to find an $h  \in \ell^{+}_{1}$ for which the potential function  $f=Uh$ satisfies the   following dichotomy for each $l\ge 0$:
\begin{equation}
 \mbox{Either 
  $F_{l}( f )\neq 0  \quad $ or  $\quad U_{I_l}$ has the form (\ref{nsz.00}).}\label{3.8mm}
 \end{equation}
 To find $h$ we take any function $h^{*}\in \ell_{1}^{+}$ and   define   successively  $h^{(n)}\in \ell^{+}_{1}$, $n\ge -1$, such that $h^{(-1)} =h^{\ast} $ and    
 \be 
  h^{(n+1)}_{j}=h^{(n)}_{j},\hspace{.2 in} \forall j\notin I_{n}, \quad\mbox{and}\quad 0<{1 \over 2} h_{j}^{\ast}\leq h^{(n)}_{j}\leq 2h_{j}^{\ast},\quad  j\ge 1\label{nsz.01}, 
\ee
and such that  $f^{(n)}:=Uh^{(n)}$ satisfies,
\be  |F_{l}( f^{(n+1)} )- F_{l}( f^{(n)} )|\leq \displaystyle\frac{|F_{l}( f^{(l+1)} ) |}{ 2^{n+2}},\quad  n\ge l+1 .\label{nsz.02}
\ee
As we point out just below (\ref{nsz.2}),  if $U_{I_l}$ is of the form (\ref{nsz.00}), (\ref{nsz.02}) is satisfied trivially   since $ F_{l}( f)=0$ for all $f$. However, when $U_{I_l}$ is not of the form (\ref{nsz.00})
we also require that $h^{(l+1)}$ is such that 
\begin{equation}
 F_{l}( f^{(l+1)} )\neq 0.\label{nsz.02m}
\end{equation}
(The actual  construction of   $\{h^{(n)};  n\ge -1\}$ is given later in this proof.)

By (\ref{nsz.01}), $\|h^{(n)} - h^{(m)}\|_{1}\leq 2\sum_{j=m}^{n} h_{j}^{*} $  for any $n>m$, hence  $h=\lim_{n\to \ff}h^{(n)}$ exists in $\ell_{1}^{+}$.
We set  $f=Uh$ and note that
\begin{equation}
   |f_{j}-f_{j}^{(n)}|=|(U(h-h^{(n)}))_{j}|\le U_{j,j}\|h-h^{(n)}\|_{1}\label{3.12mm}.
   \end{equation}
Here we use the property  pointed out in the proof of  Theorem \ref{theo-borelNS} that    
$
   U_{i,j}\le U_{i,i}\wedge U_{j,j} .\label{3.13m}
$

It follows from (\ref{3.12mm}) that
    $f_{j}=\lim_{n\to \ff}f_{j}^{(n)}$ for each $j\ge 1$ and consequently, by (\ref{nsz.02}),

\begin{equation}
    |F_{l}( f )- F_{l}( f^{(l+1)} )|\le \sum_{k=l+1}^{\ff}    |F_{l}( f ^{(k+1)})-F_{l}( f ^{(k)} )|  \leq {|F_{l}( f^{(l+1)} ) | \over 2 }.\label{3.13mm}
   \end{equation}
We see from this that  when $U_{I_l}$ is not of the form (\ref{nsz.00}),  it follows from (\ref{nsz.02m}) and (\ref{3.13mm}) that    $
 F_{l}( f   )\neq 0.$ This implies that (\ref{3.8mm}) holds.
 
 \medskip	We now describe how the $h^{(j)}$,  $j=0,1,\ldots$ are chosen. 
Assume that $h^{(-1)},\ldots, h^{(n)}$ have been chosen. We choose   $h^{(n+1)}$ as follows: 
If either $F_{n}( f^{(n)})\neq 0$ or $U\big |_{I_{n}\times I_{n}}$ has the form (\ref{nsz.00}), we set $h^{(n+1)}=h^{(n)}$.

Assume then that $F_{n}( f^{(n)} )= 0$.    
  If 
  $U _{I_{n}} $ does not have the form of (\ref{nsz.00}), it follows from the proof of Lemma \ref{lem-3.1mm}  that for all  $\ep_{p}\downarrow 0$, there exists a  
   $(g_{1,p},g_{2,p},g_{3,p})\in  B_{\ep_{p}}(f^{(n)}_{i_{n,1}},f^{(n)}_{i_{n,2}},f^{(n)}_{i_{n,3}})$ such that 
    $F_{n}( g_{1,p},g_{2,p},g_{3,p} )\ne 0$. We choose   $f^{(n+1)}=f^{(n)}$ for all indices except $i_{n,1},i_{n,2},i_{n,3}$ and $f^{(n+1)}_{ i_{n,1}},f^{(n+1)}_{ i_{n,2}},f^{(n+1)}_{ i_{n,3}}$ to be equal to one of these triples $(g_{1,p},g_{2,p},g_{3,p})$.   This gives (\ref{nsz.02m}) for $l=n$. Since  $\ep_{p}\downarrow 0$   we can take $f^{(n+1)}$ arbitrarily close to $f^{(n)}$ so that it satisfies (\ref{nsz.02}). 
    
 	   As in the proof of Theorem \ref{theo-borelNS}  we can solve the equation 
  \begin{equation}
  f^{(n+1)}_{i_{n,j}}=  f_{i_{n,j}}^{(n )}+ \sum_{k=1}^{3}U_{i_{n,j},i_{n,k}} c_{i_{n,k}},\qquad j=1,2,3. \label{3.14mmq}
   \end{equation} 
 for   $  c_{i_{n,1}},c_{i_{n,2}}, c_{i_{n,3}} $.  
      To obtain $h^{(n+1)}$  we  set  $h^{(n+1)}_{q}=h^{(n)}_{q}$ for all $q\notin I_{n}$ and for $q\in I_{n}$ we take 
 \be
 h^{(n+1)}_{q}=h^{(n)}_{q}+c^{(n)}_{q} \label{8.15mm}.
 \ee
where  $c^{(n)}_{q}$ has all its components equal to zero except for the three components $ c_{i_{n,1}},c_{i_{n,2}}, c_{i_{n,3}} $.      By taking $\ep_{p}$ sufficiently small  we can choose $ c_{i_{n,1}},c_{i_{n,2}}, c_{i_{n,3}} $ so that the third statement in   (\ref{nsz.01}) holds.  

We set $ f^{(n+1)}=U  h^{(n+1)} $ and  note that this is consistent with (\ref{3.14mmq}). \qed

   	 	 \noindent{\bf  Proof of Theorem  \ref{theo-borelNSmm} } It is clear from Theorem  \ref{theo-borelNS}   that if $(i)$ holds then $U$ is asymptoticly symmetrizable,   because in this case $\{U_{t_{i},t_{j}} \}_{i,j=1}^{k} $ is symmetrizable  for all  distinct $t_{1},\ldots ,t_{k}$ greater than or equal to $n_{0}$, for all $k$.
	 
	 Suppose that $(i)$ does not hold. Then, as in the proof of Theorem \ref{theo-borelNS}, we can find a sequence $\{n_{k};k\in \mathbb N\}$ such that $n_{k}\to \ff$ and a sequence of triples      $3n_{k}< { t_{k,1}}, 	 { t_{k,2}}, { t_{k,3}}\le 3n_{k+1}$, such that $U_{\{ t_{k,1}, t_{k,2}, t_{k,3}\}}$ does not have   all of its off diagonal elements equal.   We interchange  the   indices  ${t_{k,1}}, 	 { t_{k,2}}, { t_{k,3}}$ with the indices in  $I_{n_{k}}$; (see Lemma \ref{lem-borelNSmm}).  We can now use Lemma \ref{lem-borelNSmm} to show that $(ii)$ holds. \qed
	 
   \noindent{\bf Proof of Theorem \ref{theo-1.4} }   Let $S'=\{x_{0},x_{1},x_{2},\ldots \}$  with $\lim_{k\to \ff}x_{k}=x_{0}$. Assume that for some    integer $n_{0}$   
	 \begin{equation}
 u(  x_{j},x_{k}  )=\La_ {j}\de_{x_{j},x_{k}}+ d,\qquad \forall j,k\ge n_{0}.\label{3.28mm}
   \end{equation}
Then  $ u(  x_{j},x_{j}  ) =\La_ {j} + d$,  and since, by hypothesis,   $u(x,y)$ is continuous,  
\be
 \lim_{j\to\ff}u(  x_{j},x_{j}  ) =u(  x_{0},x_{0}  ) , \label{unif}
 \ee
 which implies that   limit $\La_ {0 }:=\lim_{j\to \ff}\La_ {j }$ must exist and  
  \begin{equation}
u(  x_{0},x_{0}  )=\La_ {0 }+d.\label{3.28gh}
  \end{equation}
  It also follows from  (\ref{3.28mm})  that $  u(  x_{j},x_{k}  ) =d$ for all $n_{0}\leq j<k$. In addition, since $\lim_{k\to\ff} u(  x_{j},x_{k}  )= u(  x_{j},x_{0}  )$,
we see that for all $j\ge n_{0} $,
     \begin{equation}
u(  x_{j},x_{0}  )=d.\label{3.28gh2}
  \end{equation}
  Comparing the last two displays we get that for all $j\ge n_{0} $,
   \begin{equation}
u(  x_{0},x_{0}  )- u(  x_{j},x_{0}  )=\La_ {0 }.\label{3.28gh3}
  \end{equation}
   This  contradicts (\ref{3.28mm}), 
because the assumption that $  u(  x_{0},x_{0}  )>u(  x_{j},x_{0}  )$
implies that $\La_ {0 }>0$, whereas  the assumption that $u$ is continuous and  (\ref{3.28gh3}) implies  that $\La_ {0 }=0$. 

Since  (\ref{3.28mm}) does not hold for any integer $n_{0}$, 
(\ref{1.9mm}) follows from  Theorem  \ref{theo-borelNSmm}.    The fact that $f$ is continuous at $x_{0}$ follows from the Dominated Convergence Theorem since    $\lim_{j,k\to\ff} u(  x_{j},x_{k}  )=u(x_{0},x_{0})$ implies that $\{u(x,y);x,y\in S'\}$ is uniformly bounded.  
    \qed

  	 \section{Proof of Lemma \ref{lem-1.1mm}  and Examples \ref{ex-1.1} and \ref{ex-1.2}}   	\label{sec-4}
 \noindent{\bf Proof of Lemma \ref{lem-1.1mm} }   
(i)   Let $m_{1},m_{2},m_{3}$ be  increasing integers such that  $m_{2}-m_{1}=m_{3}-m_{2}$ and $u(m_{2}-m_{1}) \ne u(m_{3}-m_{1}) $ and consider the  $3\times 3$ T\" oeplitz matrix
    \be
  \left (
\begin{array}{  ccc}  
 u(0) +f(m_{1})& u(m_{2}-m_{1}) +f(m_{2})  &  u(m_{3}-m_{1})  +f(m_{3}  )  \\
 u(m_{2}-m_{1}) +f(m_{1})& u(0) +f(m_{2})  & u(m_{2}-m_{1})  +f(m_{3}  )  \\
   u(m_{3}-m_{1})  +f(m_{1})&  u(m_{2}-m_{1})  +f(m_{2}) &u(0) +f(m_{3}  )   \end{array}\right )\label{pp}.
      \ee  
        By Lemma \ref{lem-1.1n}, if $\{\wt u^{f}(j,k);j,k\in \mathbb N\}$ is symmetrizable we must have
         \bea
  &&\hspace{-.3in} ( u(m_{2}-m_{1})+f(m_{2}) )(u(m_{2}-m_{1})+f(m_{3}  )  )  (    u(m_{3}-m_{1})+f(m_{1}))\label{4.2mm}\\\nn &&\quad\hspace{-.3in}=(u(m_{3}-m_{1})+f(m_{3}  ) )( u(m_{2}-m_{1})+f(m_{1}))( u(m_{2}-m_{1})+f(m_{2}) ).
   \eea
   Note  that we can cancel the term $u(m_{2}-m_{1})+f(m_{2})$ from each side of (\ref{4.2mm}) and rearrange it  to get  
   \begin{equation}
 (u(m_{2}-m_{1})- u(m_{3}-m_{1}))  (  f(m_{1})-f(m_{3})) =0.   \end{equation} 
This is not possible because   $u(m_{2}-m_{1}) \ne u(m_{3}-m_{1 }) $ and $f(m_{1})\ne f(m_{3})$.  

Since this holds for all $m_{1},m_{2},m_{3}$ satisfying the conditions above we see that  Lemma \ref{lem-1.1mm} $(i)$ holds.

 \medskip	(ii)   Consider $s_{j}\wedge s_{k}$ at the three different values,  $s_{j_{1}},s_{j_{2}},s_{j_{3}}$,  and the matrix
 \be  
  \left (
\begin{array}{  ccc}  
s_{j_{1}}+f(s_{j_{1}})&s_{j_{1}}+f(s_{j_{2}})    &s_{j_{1}}+f(s_{j_{3}}  )  \\
 s_{j_{1}}+f(s_{j_{1}})&    s_{j_{2}}+f(s_{j_{2}})  & s_{j_{2}}+f(s_{j_{3}}  )  \\
    s_{j_{1}}+f(s_{j_{1}})&  s_{j_{2}}+f(s_{j_{2}}) & s_{j_{3}}+f(s_{j_{3}}  )   \end{array}\right )\label{jj}.
      \ee   
          By Lemma \ref{lem-1.1n}, if   $ \wt v^{f}_{s_{j_{1}},s_{j_{2}},s_{j_{3}}}$ is symmetrizable we must have
 \be  ( s_{j_{1}}+f(s_{j_{2}}) )(s_{j_{2}}+f(s_{j_{3}}  )  )  ( s_{j_{1}}+f(s_{j_{1}})) =( s_{j_{1}}+f(s_{j_{3}}  ) )(  s_{j_{1}}+f(s_{j_{1}}))( s_{j_{2}}+f(s_{j_{2}}) )
   \ee   
     or, equivalently,
   \begin{equation}
 (s_{j_{1}}- s_{j_{2}})  (  f(s_{j_{3}})-f(s_{j_{2}})) =0.   \end{equation} 
 Since $s_{j_{1}}\ne s_{j_{2}}$ and  $f(s_{j_{3}})\ne f(s_{j_{2}})$ this is not possible. Therefore,    $ \wt v^{f}_{s_{j_{1}},s_{j_{2}},s_{j_{3}}}$ is not  symmetrizable.         \qed
  
  {\bf Proof of Example \ref{ex-1.1} } {\rm Let $s_{0}\in S$. We choose a sequence $s_{j}\to s_{0}$ with the property that it contains a subsequence $\{s_{j_{k}}\}$, $s_{j_{k}}\to s_{0}$, such that  
   \begin{equation}
  s_{j_{3k+1} }-  s_{j_{3k} }=  s_{j_{3k+2}}-  s_{j_{3k+1} }=a_{k}, \qquad k\ge 1.
   \end{equation}
 The kernel of 
 the $3\times 3$ matrix 
 \begin{equation}
   \wh u^{f}(s_{j_{3k +p}},s_{j_{3k +q}}),\qquad p,q=0,1,2, \label{4.8mm}
   \end{equation} 
 is
  \be
  \left (
\begin{array}{  ccc}  
1 +f(s_{j_{3k}  })& e^{-\la a_{k}} +f(s_{j_{3k+1 }})  &  e^{-\la 2a_{k}}  +f(s_{j_{3k+2} }  )  \\
e^{-\la a_{k}} +f(s_{j_{3k}  })&1 +f(s_{j_{3k+1} })  & e^{-\la a_{k}} +f(s_{j_{3k+2} }  )  \\
e^{-\la 2a_{k}} +f(s_{j_{3k}  })&  e^{-\la a_{k}} +f(s_{j_{3k+1} }) &1 +f(s_{j_{3k+2} }  )   \end{array}\right )\label{ppq},
      \ee  
   similar to (\ref{pp}).   Therefore, following the proof of Lemma \ref{lem-1.1mm}, we see that the kernel in (\ref{4.8mm}) is not symmetrizable.
 Since this holds along the subsequence  $\{s_{j_{k}}\}$, $s_{j_{k}}\to s_{0}$, we see that $ \{ \wh u_{f}(s,t); s,t\in S\}$ is not  asymptoticly symmetrizable at $s_{0} $.     
       
 The result in (\ref{1.33k}) is proved similarly. \qed

     \noindent {\bf Proof of Example \ref{ex-1.2} } The proof of Example \ref{ex-1.2} is similar to the proof of Example \ref{ex-1.2} but even simpler. This is because for all distinct values, $s_{j_{1}},s_{j_{2}},s_{j_{3}}$,    the matrix  in (\ref{jj}) is not symmetrizable.

{\footnotesize

\noindent
\begin{tabular}{lll} &  Michael Marcus
     & \hskip20pt  Jay Rosen\\  &  253 West 73rd. St., Apt. 2E
     & \hskip20pt Department of Mathematics \\ &  New York, NY 10023
  & \hskip20pt College of Staten Island, CUNY \\    &  U.S.A.     & \hskip20pt  Staten Island, NY
10314
 \\    &   mbmarcus@optonline.net 
     & \hskip20pt U.S.A. \\   &   
& \hskip20pt  jrosen30@optimum.net
\end{tabular}
\bigskip

}

  \end{document}